\documentclass[11pt,a4paper]{article}
\usepackage{graphicx,amsmath,amsfonts}
\newtheorem{theorem}{Theorem}[section]

\newtheorem{definition}[theorem]{Definition}

\newtheorem{lemma}[theorem]{Lemma}

\newcommand{\Gsub}{\underline{\Gamma}}

\def\KK{{\mathbb K}}

\newcommand{\beqa}{\begin{eqnarray}}
\newcommand{\eeqa}{\end{eqnarray}}

\newcommand{\cA}{{\cal {A}}}
\newcommand{\cC}{{\cal {C}}}
\newcommand{\cB}{{\cal {B}}}
\newcommand{\cH}{{\cal {H}}}

\newcommand{\id}{{\mathrm{id}}}

\newcommand{\ba}{\begin{array}}
\newcommand{\ea}{\end{array}}

\newcommand{\veps}{\varepsilon}

\def\>{\rangle}
\def\<{\langle}

\begin{document}

\title{COMBINATORIAL HOPF ALGEBRAS IN (NONCOMMUTATIVE) QUANTUM FIELD THEORY}

\author{ A. TANASA,\\
CPhT, CNRS, UMR 7644, \\
\'Ecole Polytechnique, 91128 Palaiseau, France, EU,
\vspace{1mm}\\ 
DFT, IFIN-HH,\\
P.O.Box MG-6, 077125 Bucure\c sti-M\u agurele, Rom\^ania, EU\\
\vspace{5mm}
E-mail: adrian.tanasa@ens-lyon.org}

\date{}
\maketitle
\begin{center}
\end{center}

\begin{abstract}
We briefly review the r\^ole played by algebraic structures like combinatorial Hopf algebras in the renormalizability of (noncommutative) quantum field theory. After sketching the commutative case, we analyze the noncommutative Grosse-Wulkenhaar model. 
\end{abstract}

{\centering\section{INTRODUCTION \label{int} }}

The r\^ole of combinatorial algebraic structures in quantum field theory has been thoroughly investigated in the last decade. Thus, the Connes-Kreimer Hopf algebra \cite{ck} is known to encode the combinatorial structure of perturbative renormalizability in commutative field theory (for a review of this topic - and of many others - the interested reader may report himself to \cite{book}). In order to investigate this algebraic structure, pre-Lie and Lie algebras can be also defined.
Going further in the understanding of quantum field theory, one can use Hochschild cohomology to guide the way from perturbative to non-perturbative physics  \cite{bk, yeats}.

When considering field theory on noncommutative Moyal space, the Grosse-Wulkenhaar model \cite{GW} was a first proposition for a renormalizable scalar model. 
Recently, a translation-invariant scalar model was proposed and also proved renormalizable on Moyal space \cite{GMRT}. 
For both these types of models, several field theoretical techniques have been recently implemented (see \cite{GW2}, \cite{io} and references within). The second type of renormalizable model has also been proved to extend also on quantum field theories based on the noncommutative Wick-Voros product \cite{io-pat}.

The combinatorial algebraic structures mentioned above have then been implemented in \cite{fab} and \cite{io-kreimer} for both these types of noncommutative models: the Grosse-Wulkenhaar one and the translation-invariant one.

\medskip

This review is structured as follows. 
The second part is a brief reminder of the definitions of the notions of algebra, bialgebra {\it etc.} used in the sequel.
In the following part, the r\^ole played by these algebraic structures  is presented for commutative  quantum field theory. In the third part, 
we present the extension of these notions in order to describe the renormalizable Grosse-Wulkenhaar model on the noncommutative Moyal space.
These last two parts follow closely \cite{fab}.
The perspective section is then devoted to the extension of these tools to the study of recent models of quantum gravity.

{\centering\section{ALGEBRAIC REMINDER}}

{\centering\subsection{ALGEBRAS}}

For further details on (Lie) algebras, the interested reader may report himself for example to \cite{teza}.

\begin{definition}[Algebra]\label{defn:algebra}
  A unital associative algebra $\cal A$ over a field $\mathbb K$ is a $\mathbb K$-linear space endowed with two algebra homomorphisms: 
  \begin{itemize}
  \item a product $m :\cA\otimes\cA\to\cA$ satisfying the \emph{associativity} condition:
    \begin{align}
      \forall\Gamma\in\cA,\,m\circ(m\otimes\id)(\Gamma)=&m\circ(\id\otimes m)(\Gamma),\label{eq:asso}
    \end{align}
  \item a unit $u :\KK\to\cA$ satisfying:
    \begin{align}
      \forall \Gamma\in\cA,\,m\circ(u\otimes\id)(\Gamma)=&\Gamma=m\circ(\id\otimes u)(\Gamma).
    \end{align}
  \end{itemize}
\end{definition}

{\centering\subsection{Hopf algebras}}

For further details on Hopf algebras, the interested reader  can refer
for example to \cite{Kassel} or \cite{Dascalescu}.

\begin{definition}[Coalgebra]\label{defn:coalgebra}
  A coalgebra $\cC$ over a field $\KK$ is a $\KK$-linear space endowed with two algebra homomorphisms: 
  \begin{itemize}
  \item a coproduct $\Delta :\cC\to\cC\otimes\cC$ satisfying the \emph{coassociativity} condition:
    \begin{align}
      \forall\Gamma\in\cC,\,(\Delta\otimes\id)\circ\Delta(\Gamma)=&(\id\otimes \Delta)\circ\Delta(\Gamma),\label{eq:coasso}
    \end{align}
  \item a counit $\veps :\cC\to\KK$ satisfying:
    \begin{align}
      \forall \Gamma\in\cC,\,(\veps\otimes\id)\circ\Delta(\Gamma)=&\Gamma=(\id\otimes\veps)\circ\Delta(\Gamma).    \end{align}
  \end{itemize}
\end{definition}

\begin{definition}[Bialgebra]
  A bialgebra $\cB$ over a field $\KK$ is a $\KK$-linear space endowed with both an algebra and a coalgebra structure (see Definitions \ref{defn:algebra} and \ref{defn:coalgebra}) such that the coproduct and the counit are unital algebra homomorphisms (or equivalently the product and unit are coalgebra homomorphisms):
  \begin{subequations}
    \label{eq:compatibility}
    \begin{align}
      \Delta\circ m_{\cB}=&m_{\cB\otimes\cB}\circ(\Delta\otimes\Delta),\ \Delta(1_\cB)=1_\cB\otimes 1_\cB,\\
      \veps\circ m_{\cB}=&m_{\KK}\circ(\veps\otimes\veps),\ \veps(1_\cB)=1.
      \end{align}
    \end{subequations}
\end{definition}

\begin{definition}[Graded Bialgebra]
  A graded bialgebra is a bialgebra graded as a linear space: 
  \begin{align}
      \cB=\bigoplus_{n=0}^\infty\cB^{(n)}   
  \end{align}
  such that the grading is compatible with the algebra and coalgebra structures:
  \begin{align}
    \cB^{(n)}\cB^{(m)} \subseteq \cB^{(n+m)}\text{ and }\Delta\cB^{(n)}\subseteq\bigoplus_{k=0}^n\cB^{(k)}\otimes\cB^{(n-k)} .
  \end{align}
\end{definition}

\begin{definition}[Connectedness]
  A connected bialgebra is a graded bialgebra $\cB$ for which $\cB^{(0)}=u(\KK)$.
\end{definition}

Let us  finally define a Hopf algebra as:

\begin{definition}[Hopf algebra]
  A Hopf algebra $\cH$ over a field $\KK$ is a bialgebra over $\KK$ equipped with an antipode map $S:\cH\to\cH$ obeying:
  \begin{align}
    m\circ(S\otimes\id)\circ\Delta=&u\circ\veps=m\circ(\id\otimes S)\circ\Delta.
  \end{align}
\end{definition}

We now end this section by recalling a lemma which will prove useful in the sequel:

\begin{lemma}[\cite{Manchon}]\label{lem:Sfree}
  Any connected graded bialgebra is a Hopf algebra whose antipode is given by $S(1_\cB)=1_\cH$ and recursively by any of the two following formulas for $\Gamma\neq 1_\cH$:
  \begin{subequations}
    \begin{align}
      S(\Gamma)=&-\Gamma-\sum_{(\Gamma)}S(\Gamma')\Gamma'',\label{eq:Srecurs}\\
      S(\Gamma)=&-\Gamma-\sum_{(\Gamma)}\Gamma'S(\Gamma'')
    \end{align}
  \end{subequations}
where we used Sweedler's notation.
\end{lemma}

{\centering\section{COMBINATORIAL HOPF ALGEBRAS OF GRAPHS - COMMUTATIVE QUANTUM FIELD THEORY}}

A graph $\Gamma$ is defined as a set of vertices and of edges (internal or external) together with an incidence relation between them.
Let us remark here that usually graph theorists and physicists have different terminologies.

\begin{figure}[bth]
\centerline{\includegraphics[width=4cm, angle=-90]{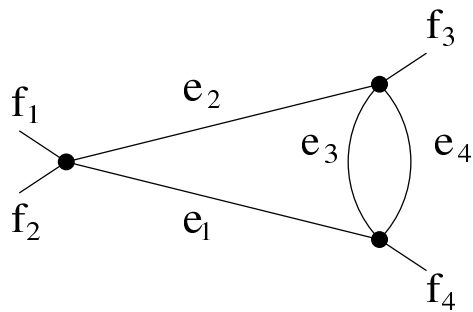}}
\caption{An example of a graph of the $\Phi^4$ quantum field theoretical model. One has $3$ vertices, $4$ internal edges ($e_1,\ldots,e_4$) and $4$ external edges ($f_1,\ldots,f_4$).
\label{ex-phi4}}
\end{figure}

In quantum field theory, one has Feynman graphs and this set of vertices and edges is given by the particle content of the model and by the type of interactions one wants to consider. For example, in the commutative $\phi^{4}_{4}$ theory, one has a unique type of particle - some scalar field $\phi$. This corresponds to a unique type of edge. The vertex type is also unique; it has valence $4$, corresponding to a $\phi^4$ interaction in quantum field theory.
An example of such a graph is the one of Fig. \ref{ex-phi4}.

\begin{definition}[Subgraph]\label{def:subgraph}
  Let $\Gamma\in\cH$, $\Gamma^{[1]}$ its set of internal edges and $\Gamma^{[0]}$ its vertices. A \textbf{subgraph} $\gamma$ of $\Gamma$, written $\gamma\subset\Gamma$, consists in a subset $\gamma^{[1]}$ of $\Gamma^{[1]}$ and the vertices of $\Gamma^{[0]}$ hooked to the edges in $\gamma^{[1]}$.
\end{definition}

Let us remark that, through such a definition, the subgraph $\gamma$ can only be truncated (that is, the external edges are not considered).

\begin{definition}
\label{def:FeynmanRules}
  The (unrenormalized) \textbf{Feynman rules} are an homomorphism $\mathbf{\phi}$ from $\cH$ to some target space $\cA$. The precise definition of $\cA$ depends on the regularization scheme employed (in the dimensional regularization scheme often used in quantum field theory, $\cA$ is the Laurent series).
\end{definition}

Let us also remark that this target algebra is naturally equipped with some Rota-Baxter algebraic structure. Recently it was studied in \cite{patras} the algebraic implications one has when relaxing this particular condition.

In quantum field theory, a special r\^ole in the process of renormalization is played by the {\bf primitively divergent graphs}. In the particular case of the $\Phi^4$ model, this class of graphs is formed by the graphs with $2$ or $4$ external edges that do not contain as subgraph any graph of $2$ or $4$ external edges.

\begin{definition}
\label{def:projectionT}
  The \textbf{projection} $\mathbf{T}$ is a map from $\cA$ to $\cA$ which has to fulfill: $\forall\Gamma\in\cH$, $\Gamma$ primitively divergent 
  \begin{align}
    (\id_{\cA}-T)\circ\phi(\Gamma)<\infty.\label{eq:conditionT}
  \end{align}
  This means that if $\phi(\Gamma)$ is divergent then its overall divergence is totally included in $T\circ\phi(\Gamma)$.
\end{definition}

We denote from now on by $\Gsub$ the set of primitively divergent subgraphs of the graph $\Gamma$.

\begin{lemma} (Lemma $3.2$ of \cite{fab})
\label{lem:coassociative1}
  Let $\Gamma\in\cH$. Provided
 \begin{enumerate}
 \item $\forall\gamma\in\Gsub,\,\forall\gamma'\in\underline{\gamma}$, $\gamma$, $\gamma'$ and $\gamma/\gamma'$ are primitively divergent,
\label{item:condition1}
 \item $\forall\gamma_{1}\in\cH,\,\forall\gamma_{2}\in\cH$ such that $\gamma_1$ and $\gamma_2$ primitively divergent, there exists gluing data $G$ such that $(\gamma_{1}\circ_{G}\gamma_{2})$ primitively divergent,\label{item:condition2}
\end{enumerate}
  the following coproduct is coassociative
\begin{subequations}
  \label{eq:coproduct1}
  \begin{align}
\Delta\Gamma=&\Gamma\otimes 1
+1\otimes\Gamma+\Delta'\Gamma,\label{eq:Deltaprime}\\
\Delta'\Gamma=&\sum_{\gamma\in\Gsub}\gamma\otimes\Gamma/\gamma.
\end{align}
\end{subequations}
\end{lemma}

Note that we denote by $\gamma/\gamma'$ the {\bf cograph} obtained by shrinking the subgraph $\gamma'$ inside $\gamma$. Shrinking a subgraph means erasing its internal structure:
\begin{itemize}
\item shrinking a subgraph with with $2$ external edges means that we replace it with a simple edge,
\item shrinking a subgraph with with $4$ external edges means that we replace it with a simple vertex.
\end{itemize}

By {\bf gluing data} we understand a bijection between the external edges of the graph to be inserted and the edges of the propagator (or vertex) where the insertion is done.

\medskip

\noindent
{\it Proof.}  Let us first remark that $(\Delta\otimes\id)\Delta=(\id\otimes\Delta)\Delta\iff (\Delta'\otimes\id)\Delta'=(\id\otimes\Delta')\Delta'$; this implies that all the following subgraphs can be considered as neither full nor empty. Let $\Gamma$ a generator of $\cH$,
 \begin{align}
    (\Delta'\otimes\id)\Delta'\Gamma=&(\Delta'\otimes\id)\sum_{\gamma\in\Gsub}\gamma\otimes\Gamma/\gamma\\
    =&\sum_{\gamma\in\Gsub}\sum_{\gamma'\in\underline{\gamma}}\gamma'\otimes\gamma/\gamma'\otimes\Gamma/\gamma\label{eq:coprodproofLeft1}\\
    (\id\otimes\Delta')\Delta'\Gamma=&\sum_{\gamma'\in\Gsub}\sum_{\gamma''\in\underline{\Gamma/\gamma'}}\gamma'\otimes\gamma''\otimes(\Gamma/\gamma')/\gamma''.\label{eq:coprodproofRight1}
  \end{align}
Using Definition \ref{def:subgraph} it is clear that $\gamma'\in\underline{\gamma}$ and $\gamma\in\Gsub$ implies $\gamma'\in\Gsub$. This implicitly uses the fact that, when shrinking a graph, what is left is independant of the surrounding of this graph and depends only on the graph itself. Equation (\ref{eq:coprodproofLeft1}) can then be rewritten as
  \begin{align}
    (\Delta'\otimes\id)\Delta'\Gamma=&\sum_{\gamma'\in\Gsub}\;\sum_{\gamma\in\Gsub}\gamma'\otimes\gamma/\gamma'\otimes\Gamma/\gamma.\label{eq:coprodproofLeft2}
  \end{align}
where in the second sum in the RHS above $\gamma'$ is considered to be a subgraph of $\gamma$ ($\gamma\ne \gamma'$). 
We just now need to prove equality between (\ref{eq:coprodproofRight1}) and (\ref{eq:coprodproofLeft2}) at fixed $\gamma'\in\Gsub$. Let us first fix a subgraph $\gamma\in\Gsub$ such that $\gamma'\subset \gamma$ ($\gamma'\ne\gamma$) and prove that there exists a graph $\gamma''\in\underline{\Gamma/\gamma'}$ such that $\gamma/\gamma'\otimes\Gamma/\gamma=\gamma''\otimes(\Gamma/\gamma')/\gamma''$. Of course the logical choice for $\gamma''$ is $\gamma/\gamma'$ because then $(\Gamma/\gamma')/(\gamma/\gamma')=\Gamma/\gamma$.

We only have to prove that $\gamma''=\gamma/\gamma'\in\underline{\Gamma/\gamma'}$. It is clear that $\gamma/\gamma'$ is a subset of internal lines of $\Gamma/\gamma'$. Then $\gamma/\gamma'\in\underline{\Gamma/\gamma'}$ if $\gamma$ primitively divergent and $\gamma'$ primitively divergent implies $(\gamma/\gamma')$ primitively divergent.

Conversely let us fix $\gamma''\in\underline{\Gamma/\gamma'}$ and prove that there exists $\gamma\in\Gsub$ containing $\gamma'$ such that $\gamma/\gamma'\otimes\Gamma/\gamma=\gamma''\otimes(\Gamma/\gamma')/\gamma''$. Let us write $\gamma'=\bigcup_{i\in I}\gamma'_{i}$ for the connected components of $\gamma'$. Some of these components led to vertices of $\gamma''$, the others to vertices of $(\Gamma/\gamma')\setminus\gamma''$. We can then define $\gamma$ as $(\gamma''\circ_{G_{I_{1}}}\bigcup_{i\in I_{1}}\gamma'_{i})\bigcup_{j\in I_{2}}\gamma'_{j}$ with $I_{1}\cup I_{2}=I$. 
This leads to the lemma. (QED)

\medskip

Let us now illustrate how this result fits the commutative $\phi^{4}$ model. In this field theory the divergent graphs have $2$ or $4$ external edges. 
The first condition of the lemma is trivial, since when one shrink a subgraph, the number of external edges is conserved.
Let us check condition \ref{item:condition2} of Lemma \ref{lem:coassociative1} for commutative $\phi^{4}$. We consider two graphs $\gamma_{1}$ and $\gamma_{2}$ with $2$ or $4$ external edges. We consider $\gamma_{0}=\gamma_{1}\circ_{G}\gamma_{2}$ for any gluing data $G$. Let $V_{i},\,I_{i}$ and $E_{i}$ the respective numbers of vertices, internal and external edges of $\gamma_{i},\,i\in\{0,1,2\}$. For all $i\in\{0,1,2\}$, we have
\begin{subequations}
  \label{eq:ConservationExtLines}
  \begin{align}
    4V_{i}=&2I_{i}+E_{i}\\
    V_{0}=&
    \begin{cases}
      V_{1}+V_{2}&\text{if $E_{2}=2$}\\
      V_{1}+V_{2}-1&\text{if $E_{2}=4$}
    \end{cases}\\
    I_{0}=&
    \begin{cases}
      I_{1}+I_{2}+1&\text{if $E_{2}=2$}\\
      I_{1}+I_{2}&\text{if $E_{2}=4$}
    \end{cases}
  \end{align}
\end{subequations}
which proves that $E=E_{1}$. Then as soon as $\gamma_{1}$ is primitively divergent so does $\gamma_{0}$. Concerning condition \ref{item:condition1} note that $\gamma''=\gamma/\gamma'\Longleftrightarrow$ it exists $G {\mbox{ s. t. }}\gamma=\gamma''\circ_{G}\gamma'$ which allows to prove that condition \ref{item:condition1} also holds and that the coproduct (\ref{eq:coproduct1}) is coassociative.

\bigskip

Consider now the unital associative algebra  $\cH$ freely generated by $1$PI Feynman graphs (including the empty set, which we denote by $1$). The product $m$ is bilinear, commutative and given by the operation of disjoint union. Let the coproduct $\Delta:\cH\to\cH\otimes\cH$ defined as
\begin{equation}
\Delta \Gamma = \Gamma \otimes 1 + 1 \otimes \Gamma + \sum_{\gamma\in \Gsub} \gamma \otimes \Gamma/\gamma, \ \forall\Gamma\in\cH.\label{eq:NCcoproduct}
\end{equation}
Furthermore let us define the counit $\varepsilon:\cH\to\mathbb K$:
\begin{equation}
\varepsilon (1) =1,\ \varepsilon (\Gamma)=0,\ \forall \Gamma\ne 1.
\end{equation}
This means that, for any non-trivial element of $\cH$, the counit returns a trivial answer and, equivalently, only for the trivial element of $\cH$ (the empty graph $1$), the result returned by the counit is non-trivial. 
Finally the antipode is given recursively by
\begin{align}
  S:\cH\to&\cH\label{eq:Antipode}\\
S(1)&=1,\nonumber\\
  \Gamma\mapsto&-\Gamma-\sum_{\gamma\in\Gsub}S(\gamma)\Gamma/\gamma.\nonumber
\end{align}

\medskip

One then has the following:

\begin{theorem} (Theorem $1$ of \cite{ck})
The quadruple $(\cH,\Delta,\veps,S)$ is a Hopf algebra.
\end{theorem}

Let us notice that $\cH$ is graded by the loop number.

\bigskip

Now let $f,g\in\text{Hom}(\cH,\cA)$ where $\cA$ is the range algebra of the projection $T$ (see above). The convolution product $\ast$ in $\text{Hom}(\cH,\cA)$ is defined by
\begin{align}
  f\ast g=&m_{\cA}\circ(f\otimes g)\circ\Delta_{\cH}.\label{eq:convolution}
\end{align}
Let $\phi$ the unrenormalized Feynman rules and $\phi_{-}\in\text{Hom}(\cH,\cA)$ the twisted antipode, defined as: $\forall\Gamma\in\cH$,
\begin{align}
  \phi_{-}(\Gamma)=&-T\big(\phi(\Gamma)+\sum_{\gamma\in\Gsub}\phi_{-}(\gamma)\ \phi(\Gamma/\gamma)\big).\label{eq:PhiMinus}
  \end{align}
The renormalized amplitude $\phi_{+}$ of a graph $\Gamma\in\cH$ is given by:
\begin{align}
  \phi_{+}(\Gamma)=&\phi_{-}\ast\phi(\Gamma).\label{eq:RenormValue}
  \end{align}

For the sake of completeness, let us recall that a combinatorial Hopf algebraic structure like the one described in this section can be also defined for more involved quantum field theories, like for example gauge theories \cite{gauge}.

Finally, let us also state that a slightly different combinatorial Hopf algebra - the {\it core Hopf algebra} - was defined in \cite{core} (and independently in \cite{thomas} for vacuum graphs).
Its coproduct does not sum on the class of {\it primitively divergent} subgraphs but on {\it all} subgraphs of the respective graph:
\begin{equation}
\Delta \Gamma = \Gamma \otimes 1 + 1 \otimes \Gamma + \sum_{\gamma\in \Gamma} \gamma \otimes \Gamma/\gamma, \ \forall\Gamma\in\cH.
\end{equation}
The definitions of the coproduct and counit follow in a straightforward manner. One can then verify that this structure is also a Hopf algebra one. The cohomology of this core Hopf algebra was also investigated (see \cite{core2} for further details).

{\centering\section{COMBINATORIAL HOPF ALGEBRAS FOR RIBBON GRAPHS - NONCOMMUTATIVE  QUANTUM FIELD THEORY}}

\begin{figure}[bth]
\centerline{\includegraphics[width=2cm, angle=-90]{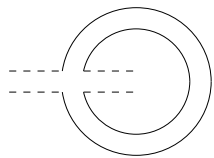}}
\caption{An example of a ribbon graph with  $1$ vertex, $1$ internal edge and $2$ external edges.
\label{ex1}}
\end{figure}

\begin{figure}[bth]
\centerline{\includegraphics[width=5cm]{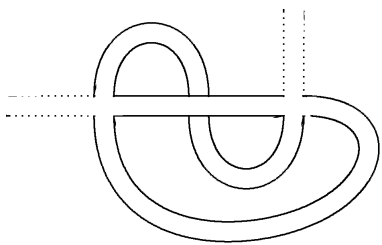}}
\caption{An example of a ribbon graph with $2$ vertices, $3$ internal edges and $2$ external edges.
\label{ex-graf-ribbon}}
\end{figure}

Feynman ribbon graph are a straightforward generalization of the Feynman graphs we have seen in the previous section. Such graphs are required in noncommutative quantum field theory. Examples of such graphs for the noncommutative $\Phi^4$ theory are given in Fig. \ref{ex1} and \ref{ex-graf-ribbon}.

More formally, a ribbon graph $\Gamma $ is an orientable surface with
boundary represented as the union of closed disks, also called
vertices, and  ribbons also called edges, such that:
\begin{itemize}
\item the disks and the ribbons intersect in disjoint line
segments,
\item each such line segment lies on the boundary of precisely one
disk and one ribbon,
\item every ribbon contains two such line segments.
\end{itemize}

One can thus think of a ribbon graph as consisting of
disks (vertices) attached to each other by thin stripes (edges) glued to
their boundaries. For any such ribbon graph $\Gamma$ there is an underlying
ordinary graph $\bar \Gamma$ obtained by collapsing the disks to points and the ribbons to edges.

Two ribbon graphs are isomorphic if there is a
homeomorphism from one to the other mapping vertices to vertices and
edges to edges. A ribbon graph is a graph with a fixed cyclic
ordering of the incident half-edges at each of its vertices.

A {\bf face} of a ribbon graph is a connected component of its boundary as a surface.
If we glue a disk along the boundary of each face we obtain a closed Riemann surface
whose genus is also called the genus of the graph.
The ribbon graph is called {\bf planar} if that Riemann surface has genus zero. The graph of Fig. \ref{ex1} is planar while the one of Fig. \ref{ex-graf-ribbon} is non-planar (it has genus $1$).

We also call a planar graph {\bf regular} if the number of faces broken by external edges is equal to $1$ and respectively {\bf irregular} if it is $>1$.

\bigskip

The definition of the Hopf algebra of non-commutative Feynman graphs which drives the combinatorics of  renormalization is formally the same as in the commutative case.  As already mentioned it has been proved that the Grosse-Wulkenhaar model is renormalizable to all orders of perturbation, so such an algebraic structure makes sense also from a ``physical'' point of view.

The class of primitively divergent graphs (one now has to sum on in the definition of the coproduct) takes explicitly into account the topology of the graph. Thus, it was proved in \cite{GW} (see also \cite{GW2} and references within) that the primitively divergent graphs of the Grosse-Wulkenhaar models are the {\it planar regular} ribbon graphs with $2$ and $4$ external edges.

\medskip

We denote by $\cH^\star$ the unital, associative algebra freely generated by $1PI$ ribbon Feynman graph. One then has the following:

\begin{theorem} (Theorem $4.1$ of \cite{fab})
\label{thm:hopf-algebra-nonComm}
The quadruple $(\cH^\star,\Delta,\veps,S)$ is a Hopf algebra.
\end{theorem}
{\it Proof.}
The only thing to prove is the coassociativity of the coproduct \eqref{eq:NCcoproduct}. Once this is done, the definition \eqref{eq:Antipode} for the antipode follows from the fact that $\cH$ is graded (by the loop number, as in the commutative case described above), connected and from Lemma \ref{lem:Sfree}.

We will use Lemma \ref{lem:coassociative1} and the fact that for all $\Gamma\in\cH$, $\Gamma$ primitively divergent is equivalent to $\Gamma$ is planar regular. Then conditions \ref{item:condition1} and \ref{item:condition2} of Lemma \ref{lem:coassociative1} are equivalent to:
\begin{enumerate}
\item for all $\gamma$ and $\gamma'\subset\gamma$ both planar regular, $\gamma/\gamma'$ is planar regular,\label{item:cond1proof}
\item for all $\gamma$ and $\gamma'\subset\gamma$ both planar regular, there exits gluing data $G$ such that $\gamma\circ_{G}\gamma'$ is planar regular.\label{item:cond2proof}
\end{enumerate}

We consider here the insertions of graphs with $4$ external edges only. 
Before proving conditions \ref{item:cond1proof} and \ref{item:cond2proof}, let us consider the insertion of a regular  graph $\gamma_{2}$ with $4$ external edges into a vertex of another graph $\gamma_{1}$.

Let $\gamma_{0}=\gamma_{1}\circ\gamma_{2}$ and for all $i\in\{0,1,2\}$ let $F_{i},I_{i},V_{i},B_{i}$ the respective numbers of faces, internal edges, vertices and broken faces of $\gamma_{i}$. 

\begin{figure}[bth]
\centerline{\includegraphics[width=2cm]{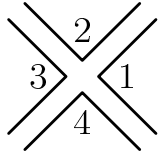}}
\caption{A ribbon vertex.
\label{fig:vert1}}
\end{figure}

\begin{figure}[bth]
\centerline{\includegraphics[width=4cm]{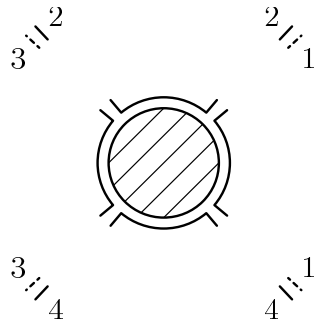}}
\caption{Insertion of a regular ribbon graph with $4$ external edges. 
\label{fig:insertion}}
\end{figure}

\begin{figure}[bth]
\centerline{\includegraphics[width=4cm]{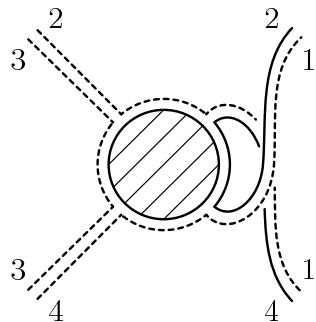}}
\caption{Another example of insertion of a ribbon graph.
\label{fig:InsRegGraph}}
\end{figure}

A ribbon vertex is drawn on figure \ref{fig:vert1}. 
We emphasize here that a cyclic ordering can be defined at the level of this vertex.
One sees that the number of faces to which the edges of that vertex belong is at most four. Some of them may indeed belong to the same face. The gluing data necessary to the insertion of $\gamma_{2}$ corresponds to a bijection between the edges of the vertex in $\gamma_{1}$ and the external edges of $\gamma_{2}$. This last one being regular ({\it i. e.} only one broken face, see above) the typical situation is represented on figure \ref{fig:insertion}. It should be clear that $F=F_{2}-1+F_{1}-n$ for some $n\geq 0$. Note that $F_{2}-1$ is the number of internal faces of $\gamma_{2}$ i.e. the number of faces of the blob. The number $n$ depends on the gluing data. It vanishes if the insertion respects the \emph{cyclic ordering} of the vertex. For example the following bijection $\sigma$ does:
\begin{align}
  \sigma((1',2'))=&(2,3),\ \sigma((2',3'))=(3,4),\ \sigma((3',4'))=(4,1),\ \sigma((4',1'))=(1,2).  
\end{align}
As in equations \eqref{eq:ConservationExtLines}, $I_{0}=I_{1}+I_{2}$ and $V_{0}=V_{1}+V_{2}-1$. It follows that the genus of $\gamma_{0}$ satisfies
\begin{align}
  g(\gamma_{0})=&g(\gamma_{1})+g(\gamma_{2})+n.\label{eq:GenusInsert}
\end{align}
Moreover by exhausting the $4!/4$ possible insertions, one checks that $B_{0}\geq B_{1}$. 
Figure \ref{fig:InsRegGraph} shows an insertion of a regular graph with $4$ external edges which increases the number of broken faces by one: now trajectories $(1,4),(1,2)$ and $3$ are external faces (edge $(2,4)$ is still an internal one).

Let us now turn to proving that the algebra of non-commutative Feynman graphs described above fulfills conditions \ref{item:cond1proof} and \ref{item:cond2proof}.
\begin{enumerate}
\item $\gamma,\gamma'$ planar implies $\gamma/\gamma'$ planar thanks to equation \eqref{eq:GenusInsert}. Furthermore $B(\gamma)=1$ implies $B(\gamma/\gamma')=1$ due to the preceding remark.
\item For condition \ref{item:cond2proof} one chooses gluing data $G$ respecting the cyclic ordering of the vertex. Then one has $g(\gamma\circ_{G}\gamma')=g(\gamma)+g(\gamma')=0$. The cyclic ordering of the insertion ensures $B(\gamma\circ_{G}\gamma')=B(\gamma)=1$.
\end{enumerate}
The proof is thus completed. (QED)

The renormalized Feynman amplitudes is defined analogously as in the commutative case (see previous section).

\medskip

Let us end this section by recalling that in \cite{io-kreimer} it was also studied in detail the r\^ole played in noncommutative quantum field theory (both for the Grosse-Wulkenhaar and for the translation-invariant model \cite{GMRT}) by the Hochschild cohomology of these types of Hopf algebras; non-trivial examples of $1$ and $2$ loops have been explicitly worked out.

{\centering\section{PERSPECTIVES - COMBINATORIAL HOPF ALGEBRAS FOR TENSOR MODELS}}

The main perspective of the combinatorial algebraic approaches presented here is their extension for the study of the renormalizability properties of quantum gravity models. 
The group field theory formalism of quantum gravity (for a general review,
see \cite{gft}) is the most adapted one for such a study, since it is formulated as a (quantum) field theory. These models were developed as a generalization of
$2-$dimensional matrix models (which naturally make use of ribbon graphs just as noncommutative quantum field theories do) to the $3-$ and $4-$dimensional cases.

The natural candidates for generalizations of matrix models in higher dimensions ($>2$)
are {\it tensor models}. In the combinatorial simplest case, the elementary cells that, by gluing
together form the space itself, are the $D-$simplices ($D$ being the dimension of space). Since a
$D-$simplex has $(D+1)$ facets on its boundary, the backbone of group field theoretical models in
$D-$dimension should be some abstract $\phi^{(D+1)}$ interaction on rank $D$ tensor fields $\phi$.

As the size of the matrix in matrix models increases, the scaling of observables
promotes the particular class of planar graphs which corresponds to triangulations of the
sphere, higher genus surfaces being exponentially suppressed.

In a recent paper 
\cite{freidel} it was
defined a contraction
procedure on a particular family of Feynman graphs of the $3$-dimensional tensor model, called Feynman
graphs of type $1$.
In commutative field theory, all
connected graphs can be contracted to points. In a noncommutative quantum field theory on
Moyal space, the planar graphs can be contracted to a non-local Moyal vertex; this does not
remain true for non-planar graphs, but renormalizability requires this only for planar graphs.
What 
it was thus proved in \cite{freidel}
is that, for the Feynman graphs of type $1$, the generalization of these contraction
procedures can be defined at the level of the $3-$dimensional tensor model

It is thus interesting to investigate weather or not these
Feynman graphs of type $1$ can be the primitively divergent Feynman graphs of  the $3-$dimensional tensor model. 
One proposition for achieving this task
is to use using the powerful tool of algebraic combinatorics that were sketched in this review. 
As explained above, in order to have an appropriate Connes-Kreimer combinatorial Hopf algebraic
structure (which, in commutative and noncommutative theories on Moyal space -
underlies renormalizability) one needs to find a way of inserting primitively divergent
Feynman graphs into primitively divergent Feynman graphs such that the resulting graph is
also primitively divergent. 

It is thus interesting to generalize this for the tensor Feynman graphs of this type of model.
Note that this is highly non-trivial because of the complexity of the topological properties of
this type of graphs. One needs to explicitly consider all pairings between the external edges of
the graph to be inserted and the edges at the insertion place. For each such possibility, one then
has to investigate if the resulting graph remains of type $1$. Moreover, when performing this
type of analysis one has to consider the symmetries of the vertex (for example, in the case of
the ribbon graphs of noncommutative field theories on Moyal space, one has a cyclical
symmetry at the vertex level, as explained in the previous section). These symmetries are responsible for identifying a priori
distinct gluing data.

One further property that one can prove using this method, is that the type 1 Feynman
graphs are the only Feynman graphs which can be inserted in such a manner.

Finally, let us also stress on the fact that different topological and analytical insights on this type of formulation of quantum gravity models have been investigated in the recent literature (see \cite{ultimele} and references within).

For the sake of completeness, 
let us also recall that, in a different type of quantum gravity formulation (the spin-foam one, see for example \cite{sf}), this type of algebraic combinatorial study has already been performed \cite{io-sf}.

{\centering\section{REFERENCES}}

{\centering\subsubsection*{Acknowledgements}}
This work was partly supported by 
the CNCSIS grant Idei 454/2009, ID- 44 and by the grant PN 09 37 01
02.

\begin{center}

\end{center}
\end{document}